\documentclass[preprint, a4paper, 11pt, oneside, onecolumn]{elsarticle}
\usepackage{amssymb, amsthm}
\usepackage{graphicx}
\usepackage{subfigure}
\usepackage{float}    
\usepackage{amsthm}
\usepackage{amsmath}
\usepackage{mathrsfs}
\usepackage[misc]{ifsym}
\usepackage{diagbox}
\usepackage{graphicx} 
\usepackage{amsmath,latexsym,amssymb,amsfonts,amsbsy}
\usepackage{subfigure} 
\usepackage{algorithm}
\usepackage{algpseudocode}
\usepackage{color,xcolor}
\usepackage{booktabs}
\usepackage{threeparttable}
\usepackage{multicol}
\usepackage{multirow}
\usepackage{epstopdf}
\numberwithin{equation}{section}
\newtheorem{theorem}{Theorem}[section]
\newtheorem{lemma}{Lemma}[section]

\newtheorem{remark}{Remark}[section]

\newtheorem{definition}{Definition}[section]

\newtheorem{proposition}[theorem]{Proposition}

\begin{document} 
\begin{sloppypar}
    \begin{frontmatter}
       \title{Numerical Analysis of Locally Adaptive Penalty Methods For The Navier-Stokes Equations}
        \author{Rui Fang}
\address{Department of Mathematics, University of Pittsburgh, Pittsburgh, PA 15260, USA}
\ead{ruf10@pitt.edu}
        \begin{abstract}
        	Penalty methods relax the incompressibility condition and uncouple velocity and pressure. Experience with them indicates that the velocity error is sensitive to the choice of penalty parameter $\epsilon$. So far, there is no effective \'a prior formula for $\epsilon$. Recently, Xie developed an adaptive penalty scheme for the Stokes problem that picks the penalty parameter $\epsilon$ self-adaptively element by element small where $\nabla \cdot u^h$ is large. Her numerical tests gave accurate fluid predictions. The next natural step, developed here, is to extend the algorithm with supporting analysis to the non-linear, time-dependent incompressible Navier-Stokes equations.  In this report, we prove its unconditional stability, control of $\|\nabla \cdot u^h\|$, and provide error estimates. We confirm the predicted convergence rates with numerical tests.
        \end{abstract}

        \begin{keyword}Navier-Stokes equations, penalty methods, adaptive algorithms, numerical analysis, FEM\\
        \vspace{0.2cm}
        \end{keyword}   
    \end{frontmatter}

\section{Introduction}\label{intro}
Consider the incompressible Navier-Stokes equations (NSE) with no-slip \\
boundary conditions:
\begin{equation}
\begin{split}
    u_t+u\cdot \nabla u+\nabla p-\nu\Delta u=f, \text{ and } \nabla \cdot u=0 \text{ in } \Omega \times [0,T]. \\
\end{split}\label{NSE}%
\end{equation}
Here, $u$ represents the velocity, $p$ the pressure, $f$ the body force, and $\nu$ the kinetic viscosity. The incompressibility condition, $\nabla \cdot u =0$, couples the velocity and pressure, making the system of equations (\ref{NSE}) more complex and challenging to solve. Penalty methods relax this condition, uncoupling the pressure and velocity and resulting in a system that is easier to solve. The penalized NSE, Temam \cite{temam1968methode}, is
\begin{equation}\label{penalty1}
\begin{gathered}
u_{\epsilon,t} -\nu \Delta u_{\epsilon} + u_{\epsilon}\cdot \nabla u_{\epsilon} + \frac{1}{2} (\nabla \cdot u_{\epsilon}) u_{\epsilon} + \nabla p_{\epsilon} =f,\\
 \nabla \cdot u_\epsilon +\epsilon p_\epsilon=0, \text{ where } 0<\epsilon \ll 1. 
\end{gathered}
\end{equation}
We can eliminate the pressure in the first equation by setting $p_\epsilon = -\frac{1}{\epsilon} \nabla \cdot u_\epsilon$ to speed up the calculation,  Heinrich and Vionnet \cite{heinrich1995penalty}.

The velocity error is sensitive to the choice of penalty parameter $\epsilon$, \\
Bercovier and Engelman 1979 \cite{bercovier1979finite}. Where $\epsilon$ is too big the system will poorly approximate incompressible flow when $\epsilon$ is too small ill-conditioning results, Hughes, Liu and Brooks \cite{hughes1979finite}. Xie \cite{xie2022adaptive} resolved this dilemma for the Stokes problem by developing a self-adaptive $\epsilon$ penalty scheme. This report studies its natural extension to the non-linear, time-dependent NSE.

The main idea behind adaptive $\epsilon$ penalty methods is to monitor $\nabla \cdot u^h_\epsilon$ and adjust the penalty parameter $\epsilon$ locally in both space and time to ensure that $\|\nabla \cdot u^h_\epsilon\|< \text{Tolerance}$, user-supplied. 
We localize the tolerance to an element-wise condition. We then monitor $\nabla \cdot u^h_\epsilon$ on each mesh element, denoted as $\|\nabla \cdot u^h_\epsilon \|_{\Delta}$, and adapt the penalty parameter for each mesh element $\epsilon_{\Delta}$. We describe the locally adaptive penalty algorithm in detail in Section \ref{sec: algorithm}.

In Section 3, we prove that the method (equation (\ref{local_penalty})) is energy stable and controls $\|\nabla \cdot u\|$ effectively. Then, we prove an error estimate for the semi-discrete approximation. In Section \ref{sec: tests}, we test these theoretical results using a known exact solution, a modified Green-Taylor vortex. Additionally, we test the method on a more complex flow scenario of fluid flows between offset cylinders. We combine the locally adaptive $\epsilon$ penalty method with adaptive time stepping and test it with a flow problem with sharp transition regions.

\subsection{Formulation}\label{formulation}
Let $\Omega$ be a bounded, open polyhedral domain in $\mathbb{R}^d$, where $d=2$ or $3$. Define the fluid velocity space and pressure space, respectively:
\begin{equation*}
    \begin{gathered}
    X:=H^{1}_{0}(\Omega)=\{ v\in L^2(\Omega)^d:\nabla v\in L^2(\Omega)^{d\times d} \text{ and } v=0 \text{ on } \partial \Omega\}, \text{ and}\\
    Q:=L^2_0(\Omega)=\{q\in L^2(\Omega): \int_{\Omega} q \, dx=0\}.
    \end{gathered}
\end{equation*}
Let $X^h \subset X$ be a finite element space for the fluid velocity, and $Q^h\subset Q$ be a finite element space for the fluid pressure. $(\cdot,\cdot)$ is the $L^2$ inner product with the norm $\|\cdot\|$. $\Delta$ denotes a mesh element, i.e. take any $\phi$, we have $\int_{\Omega}\phi\, dx =\sum_{\Delta} \int_{\Delta}\phi \, dx$, and the volume of a domain $\Omega$ is denoted as $|\Omega|$.\\
The weak formulation of the NSE is: find $u(t):[0,T]\to X$, and $p(t):[0,T] \to Q$ such that
\begin{equation}
    (u_t,v)+(u\cdot \nabla u,v)-(p,\nabla \cdot v)+\nu(\nabla u,\nabla v)=(f,v).  \label{weakNSE}
\end{equation}
The finite element approximation of the penalized NSE is to find $u^h_\epsilon:[0, T] \to X^h$, $\forall v^h\in X^h$ such that
\begin{equation}\label{local_penalty}
\begin{gathered}
(u^h_{\epsilon,t},v^h)+(u^h_\epsilon \cdot \nabla u^h_\epsilon, v^h)+ \frac{1}{2}((\nabla \cdot u^h_\epsilon)) u^h_\epsilon, v^h) +\nu (\nabla u^h_\epsilon,\nabla v^h)\\
+\sum_{\Delta} \epsilon_{\Delta}^{-1}\int_{\Delta}  \nabla \cdot u^h_\epsilon \nabla \cdot v^h\,dx=(f,v^h).
\end{gathered}
\end{equation}
Standard energy estimates suggest that halving $\epsilon$ will halve $\|\nabla \cdot u^h_\epsilon\|$ Xie \cite{xie2022adaptive}, Shen \cite{shen1995error} and He \cite{he2005optimal}. The basic idea in Xie \cite{xie2022adaptive} is to convert the global condition $\|\nabla \cdot u^h_\epsilon\|\leq TOL$ to a local one $\|\nabla \cdot u^h_\epsilon\|^2_{\Delta} \leq LocTol_{\Delta}$. Then select a smaller local $\epsilon$ for the next timestep where $\|\nabla \cdot u^h_\epsilon\|^2_{\Delta} > LocTol_{\Delta}$ and a larger $\epsilon$ where $\|\nabla \cdot u^h_\epsilon\|^2_{\Delta}<LocTol_{\Delta}$. We update $\epsilon$ by the following:
\begin{equation*}
\epsilon_{\Delta, \text{new} }= \rho\epsilon_{\Delta, \text{old}},
\end{equation*}
where 
\begin{equation*}
\rho = \frac{LocTol_{\Delta}}{\|\nabla \cdot u^h_\epsilon\|^2_{\Delta}}.
\end{equation*}
\subsection{Algorithms}\label{sec: algorithm} In Algorithm \ref{alg: elementwise}, we present the pseudo-code of the locally adaptive penalty method. Xie \cite{xie2022adaptive} gives a numerical test on the locally adaptive penalty method for the non-linear NSE (see p. 11, Algorithm 2 and test 4 in p. 15 of \cite{xie2022adaptive}). In this report, we analyze the non-linear and time-dependent NSE in detail. We prove its unconditional stability and error estimates and conduct a suite of numerical tests.
\begin{algorithm}
\caption{Locally adaptive $\epsilon$ penalty algorithm}\label{alg: elementwise}
\begin{algorithmic}
\State Given $TOL$, $\Delta t$, $EpsMin$, $EpsMax$. Set initial $\epsilon_{\Delta} =1$. Compute the local tolerance on each triangle 
$LocTol_{\Delta}=\frac{1}{2}\frac{TOL^2}{|\Omega|}|\Delta|$.
\While{$t < T$}:
\State Update $t_{n+1}=t_n+\Delta t$.
\State Given $u^h_{\epsilon,n}$, solve for $u^h_{\epsilon,n+1}$ with penalty method: find $u^h_{\epsilon,n+1} \in X^h$ such that
\begin{equation*}
\begin{gathered}
    (\frac{u^h_{\epsilon,n+1}-u^h_{\epsilon,n}}{\Delta t},v^h)+b^*(u^h_{\epsilon,n},u^h_{\epsilon,n+1},v^h)+\nu (\nabla u^h_{\epsilon,n+1},\nabla v^h)\\
    +\sum_{\Delta} \epsilon_{\Delta}^{-1}\int_{\Delta}  \nabla \cdot u^h_{\epsilon,n+1} \nabla \cdot v^h\, dx=(f^{n+1},v^h).
    \end{gathered}
\end{equation*}
\State Compute divergence on each triangle 
\begin{equation*}
est_{\Delta}=\int_{\Delta} |\nabla \cdot u^{h}_{\epsilon,n+1}|^2\, dx.
\end{equation*}
\State Update $\epsilon_\Delta$:
\begin{equation*}
\begin{gathered}
\rho = \frac{LocTol_\Delta}{est_\Delta},\\
\epsilon_{\Delta,n+2} \leftarrow min (max(EpsMin, \rho\times \epsilon_{\Delta, n+1}), EpsMax).
\end{gathered}
\end{equation*}
\State Recover pressure $p$ if needed
\begin{equation*}
p_{\Delta, n+1} = -\frac{1}{\epsilon_{\Delta,n+1}} \nabla \cdot u^{h}_{\epsilon,n+1}.
\end{equation*}
\EndWhile
\end{algorithmic}
\end{algorithm} 
\subsection{Related work} The penalty method was first introduced by Courant \cite{courant1943variational} and applied to unsteady NSE by Temam \cite{temam1968methode}. Shen \cite{shen1995error} proved error estimates for continuous time with the constant $\epsilon$ and for the backward Euler time discretization. This analysis suggested $\epsilon = \Delta t$, the timestep. Other proposals include $\epsilon= (\text{machine } \epsilon)^{1/2}$, Heinrich and Vionnet \cite{heinrich1995penalty}. He \cite{he2005optimal} and He and Li \cite{he2010penalty} gave optimal error estimates with conditions on $\epsilon, \Delta t$, and mesh size $h$. 
Recently, Layton and Xu \cite{layton2023conditioning} analyzed the conditioning of the linear system from the penalty method. Bernardi, Girault, and Hecht \cite{bernardi2003posteriori} studied the penalty method on adaptive meshes and locally adaptive $\epsilon$ in space by a different algorithm for the Stokes problem. They found adapting $\epsilon$ by their method inefficient. Some reports have studied how to pick the estimator to control the residual $\|\nabla \cdot u^h_\epsilon\|$, including Layton and Xu \cite{layton2023conditioning} $\frac{\|\nabla \cdot u\|}{\|u\|}$, Kean, Xie, and Xu \cite{kean2023doubly} $\frac{\|\nabla \cdot u\|}{\|\nabla u\|}$, and Xie \cite{xie2022adaptive} $\|\nabla\cdot u\|$ for the Stokes problem. The choice of penalty parameter is also studied in grad-div stabilization to improve the pressure-robustness, John, Linke, Merdon, Neilan, and Rebholz \cite{john2017divergence}. 
\section{Notation and preliminaries}\label{rec: preliminary} The Sobolev spaces and Lebesgue spaces on $\Omega$ are denoted by $W^{k,p}(\Omega)$ and $L^p(\Omega)$ respectively, equipped with the norms $\|\cdot\|_{k,p}$ and $\|\cdot\|_{L^p}$. We denote the semi norm on $W^{k,p}(\Omega)$ as $|\cdot|_{k,p}$. For $p=2$, we write $|\cdot|_{k}$ for short. Given a Banach space $X$, we define the norm on $L^p(0, T; X)$.
\begin{equation*}
\|\cdot\|_{L^P(X)} := \left(\int_{0}^T\|\cdot\|^p_{X}\, dt\right)^{1/p} \text{ if }1\leq p<\infty, \ \text{ess} \sup_{t\in [0,T]}\|\cdot\|_X \text{ if } p=\infty.
\end{equation*}
We summarize standard properties of finite element space, John \cite{john2016finite}, Layton \cite{layton2008introduction}, we assume hold. Assume $X^h$ satisfies the following approximation properties for $0\leq s\leq m$:
\begin{equation}\label{approximation-properties}
    \begin{split}
    \inf_{v\in X^h}\|u-v\|\leq Ch^{s+1} |u|_{s+1}, \forall u\in X \cap H^{s+1}(\Omega)^d,\\
    \inf_{v\in X^h} \|\nabla (u-v)\|\leq Ch^s|u|_{s+1}, \forall u\in X \cap H^{s+1}(\Omega)^d.
    \end{split}
\end{equation}
The space $H^{-1}(\Omega)$ denotes the dual space of all the bounded linear functionals on $X:=H^{1}_0$ equipped with the norm:
\begin{equation*}
    \|f\|_{-1} =\sup_{0\neq v \in X} \frac{(f,v)}{\|\nabla v\|}.
\end{equation*}
\begin{lemma}\label{h-inverse-projection} Let $P^h$ be the $L^2$ projection. Then
\begin{equation*}
\|u - P^hu\|_{-1} \leq C h\|u -P^hu\|.
\end{equation*}
\begin{proof}
Since $P^h$ is the $L^2$ projection, $\forall v^h \in X^h$, $(u-P^h u, v^h)= 0$. We have
\begin{equation*}
\begin{gathered}
\|u-P^h u\|_{-1} = \sup_{0\neq v\in X} \frac{(u -P^h u, v)}{ \|\nabla v\|}\\
= \sup_{0\neq v\in X}  \inf_{v^h \in X^h}\frac{(u -P^h u, v-v^h)}{ \|\nabla v\|} \leq \|u-P^h u\|\sup_{0\neq v\in X} \inf_{v^h \in X^h}\frac{\|v-v^h\|}{\|\nabla v\|}.
\end{gathered}
\end{equation*}
By the interpolation estimate. We have, using equation (\ref{approximation-properties}) with $s=0$,
\begin{equation*}
\inf_{v^h\in X^h}\| v-v^h\|\leq Ch \|\nabla v\|.
\end{equation*}
Hence $\|u - P^h u\|_{-1}\leq Ch \|u -P^h u\|$.
\end{proof}
\end{lemma}
Denote the skew-symmetric trilinear form: $\forall u,v,w \in X$,
\begin{equation*}
\begin{split}
b^*(u,v,w) :=\frac{1}{2}(u\cdot \nabla v, w)-\frac{1}{2} (u\cdot \nabla w, v). 
\end{split}
\end{equation*}
\begin{lemma}(\text{See Layton \cite{layton2008introduction} p.123 p.155}) $\forall u,v,w \in X$, the trilinear term\\
$b^*(u,v,w)$ is equivalent to
\begin{equation*}
\begin{split}
b^*(u,v,w)=(u\cdot \nabla v, w)+\frac{1}{2} \left((\nabla \cdot u) v, w\right).
\end{split}
\end{equation*}
\end{lemma}
\begin{lemma}\label{upperbounds}(Continuity and a sharper bound on the the trilinear form see Layton \cite{layton2008introduction} p.11) Let $\Omega \in \mathbb{R}^d$, where $d=2$ or $3$. $\forall u,v,w \in X$,
\begin{equation*}
\begin{gathered}
|b^*(u,v,w)|\leq C(\Omega)\|\nabla u\|\|\nabla v\| \|\nabla w\|,\\
|b^*(u,v,w)|\leq C(\Omega) \sqrt{\|u\|\|\nabla u\|}\|\nabla v\|\|\nabla w\|.
\end{gathered}
\end{equation*}
\end{lemma}

\begin{lemma}
(H\"older's and Young's inequality\label{holderandyoung})
For any $\sigma>0$, $1\leq p\leq \infty$, $\frac{1}{p}+\frac{1}{q}=1$, the H\"older and Young's inequality are as following:
\begin{equation*}
    (u,v)\leq \|u\|_{L^p}\|v\|_{L^q},\ \ \text{and }
    (u,v)\leq \frac{\sigma}{p} \|u\|_{L^p}^p+
    \frac{\sigma^{-q/p}}{q} \|v\|_{L^q}^q.
\end{equation*}
\end{lemma}
For $\Delta$ a mesh element (e.g. a triangle in 2d), let $(u,v)_{\Delta} = \int_{\Delta} u \cdot v \, dx$ and $\|\cdot\|_{L^q(\Delta)}$ denote the $L^2(\Delta)$ inner product and $L^q(\Delta)$ norm.
\begin{lemma}
(Xie \cite{xie2022adaptive} Lemma 2.2 p.4, Strong monotonicity and local Lipschitz continuity) \label{strongMonotonicity} Let $u,v,w \in X$, on each mesh element $\Delta$, there exist constants $C_1$ and $C_2$ such that the following inequalities hold respectively:
\begin{equation}
   (|\nabla \cdot u|^2 \nabla \cdot u - |\nabla \cdot w|^2 \nabla \cdot w, \nabla \cdot (u-w))_{\Delta}\geq C_1 \|\nabla \cdot (u-w)\|^4_{L^4(\Delta)}, \label{StongMonoton}
\end{equation}
\begin{equation}
    (|\nabla \cdot u|^2 \nabla \cdot u - |\nabla \cdot w|^2 \nabla \cdot w, \nabla \cdot v)_{\Delta} \leq C_2 r^2 \|\nabla \cdot (u-w)\|_{L^4(\Delta)} \|\nabla \cdot v\|_{L^4(\Delta)}, \label{Lipschitz}
\end{equation}
where $r=\max\big\{\|\nabla \cdot u\|_{L^4 (\Delta)},\|\nabla \cdot w\|_{L^4(\Delta)}\big\}$.
\end{lemma}
\begin{definition} (see Xie \cite{xie2022adaptive} p.3)
The local tolerance on an element $\Delta$, $LocTol_\Delta$, is defined as follows.
\begin{equation*}
\begin{split}
     LocTol_\Delta:= \frac{1}{2} \frac{TOL^2}{|\Omega|}|\Delta|.
     \end{split}
\end{equation*}
\end{definition}
 \begin{lemma} (See Xie \cite{xie2022adaptive} p.3)
 \label{tolerance}
If $\|\nabla \cdot u^h_\epsilon\|^2_\Delta \leq  LocTol_\Delta$, there yields
\begin{equation*}
\|\nabla \cdot u^h_\epsilon\|^2 = \sum_{\Delta} \int_{\Delta} |\nabla \cdot u^h_\epsilon|^2\, dx \leq \sum_{\Delta} LocTol_\Delta = \frac{1}{2} TOL^2.
\end{equation*}
\end{lemma}
\section{Convergence of the semi-discrete approximation}\label{sec: analysis} In this section, we derive stability bounds and error estimates for the locally adaptive penalized NSE. We assume the penalty parameter is pointwise, which means that for each mesh element $\Delta$, the penalty parameter can vary spatially across the element. We denote it as $\epsilon_{\Delta}(x)$. Since a new $\epsilon_{\Delta}$ selected each timestep, $\epsilon_{\Delta}$ will also depend on time. To simplify non-essential notation we suppose notational time dependence as $\epsilon_{\Delta}$. Specifically at a fixed timestep, for a given mesh element $\Delta$, we update the penalty parameter as follows:
\begin{equation*}
\epsilon_{\Delta, \text{new}}= \rho(x)\epsilon_{\Delta, \text{old}},
\end{equation*}
where 
\begin{equation*}
\rho(x) = \frac{LocTol_{\Delta}}{|\nabla \cdot u^h_\epsilon|^2}.
\end{equation*}
We start with $\epsilon =1$. Therefore, the first step
\begin{equation*}
\epsilon_{\Delta}(x)=\frac{LocTol_{\Delta}}{|\nabla \cdot u^h_\epsilon(x)|^2}.
\end{equation*}
Hence
\begin{equation*}
\begin{split}
p^h_\epsilon = -\frac{1}{\epsilon} \nabla \cdot u^\epsilon=-\frac{|\nabla \cdot u^h_\epsilon(x)|^2}{LocTol_{\Delta}} \nabla \cdot u^h_\epsilon.
\end{split}
\end{equation*}
It yields
\begin{equation}\label{penaltyNSE}
\begin{gathered}
     ({u^h_\epsilon}_t,v^h)+b^*(u^h_\epsilon,u^h_\epsilon,v^h)+\nu (\nabla u^h_\epsilon,\nabla v^h)\\
     +\sum_{\Delta } LocTol_{\Delta}^{-1} (|\nabla \cdot u^h_\epsilon |^2 \nabla \cdot u^h_\epsilon \nabla \cdot v^h)=(f,v^h).
     \end{gathered}
\end{equation}
We first present the stability analysis of the semi-discrete approximation as follows.
\begin{theorem}\label{stability-bound}
Let $\mathcal{T}^h$ be a mesh of $\Omega$ and $\Delta$ denote a mesh element in $\mathcal{T}^h$, then the solution to (\ref{penaltyNSE}) is stable, and the following stability bound holds:
\begin{equation}\label{ppstability}
\begin{gathered}
\|u^h_\epsilon(t)\|^2 +2\int_{0}^t \left(\sum_{\Delta} \frac{1}{LocTol_{\Delta}} \int_{\Delta} |\nabla \cdot u^h_\epsilon(t')|^4\, dx \right) \, dt' \\
+\nu\int_{0}^t \|\nabla u^h_\epsilon\|^2\, dt'\leq\|u_0\|^2+\frac{1}{\nu}\int_{0}^t \|f(t')\|^2_{-1}\,dt'. 
\end{gathered}
\end{equation}

\begin{proof}
Take $v^h=u^h_\epsilon$ in (\ref{penaltyNSE}). We have 
\begin{equation*}
    \begin{split}
     \frac{1}{2}\frac{d}{dt}\|u^h_\epsilon\|^2 +\nu\|\nabla u^h_\epsilon\|^2+\sum_{\Delta} \frac{1}{LocTol_{\Delta}} \int_{\Delta} |\nabla \cdot u^h_\epsilon|^4\, dx&=(f,u^h_\epsilon).
     \end{split}
\end{equation*}
$(f, u^h_\epsilon)\leq \|f\|_{-1}\|\nabla v^h_\epsilon\|$. Apply Lemma \ref{holderandyoung} H\"older and Young's inequality:
\begin{equation*}
\begin{gathered}
     \frac{1}{2}\frac{d}{dt}\|u^h_\epsilon\|^2 +\nu\|\nabla u^h_\epsilon\|^2+\sum_{\Delta} \frac{1}{LocTol_{\Delta}} \int_{\Delta} |\nabla \cdot u^h_\epsilon|^4\, dx\\
     \leq \frac{1}{2\nu}\|f\|^2_{-1}+\frac{\nu}{2}\|\nabla u^h_\epsilon\|^2.
     \end{gathered}
\end{equation*}
Combine similar terms. We have
\begin{equation*}
     \frac{1}{2}\frac{d}{dt}\|u^h_\epsilon\|^2 +\frac{\nu}{2}\|\nabla u^h_\epsilon\|^2+\sum_{\Delta} \frac{1}{LocTol_{\Delta}} \int_{\Delta} |\nabla \cdot u^h_\epsilon|^4\, dx\leq \frac{1}{2\nu}\|f\|^2_{-1}.
\end{equation*}
Integrate over time from $0$ to $t$, and the stability bound (\ref{ppstability}) follows.
\end{proof}
\end{theorem}
Proposition \ref{l4-div-bound} proves the global tolerance can effectively control $\nabla \cdot u^h_\epsilon$ for the locally adaptive $\epsilon$ algorithm. 
\begin{proposition}\label{l4-div-bound}
Let $TOL$ denote the global tolerance, then the solution $u^h_\epsilon$ to (\ref{penaltyNSE}) satisfies
\begin{equation}
   \int_{0}^t \|\nabla \cdot u^h_\epsilon(t')\|^4_{L^4}\, dt'\leq \frac{TOL^2 \max |\Delta|}{4|\Omega|}\left( \|u_0\|^2+ \frac{1}{\nu} \int_{0}^T\|f(t')\|^2_{-1}\, dt'\right).
\end{equation}
\begin{proof}
From Theorem \ref{stability-bound}, there holds
\begin{equation}
2\int_{0}^t \left(\sum_{\Delta} \frac{\int_{\Delta} |\nabla \cdot u^h_\epsilon(t')|^4\, dx}{LocTol_{\Delta}}  \right) \, dt' 
\leq\|u_0\|^2+\frac{1}{\nu}\int_{0}^t \|f(t')\|^2_{-1}\,dt'.
\end{equation}
Recall $LocTol_\Delta=\frac{1}{2} TOL^2 \frac{|\Delta|}{|\Omega|}$. We have
\begin{equation*}
\begin{gathered}
2\int_{0}^t \left( \sum_{\Delta} \frac{2|\Omega|}{TOL^2|\Delta|}\int_{\Delta} |\nabla \cdot u^h_\epsilon(t')|^4\, dx \right)\, dt' \leq \|u_0\|^2+\frac{1}{\nu} \int_{0}^t \|f(t')\|^2_{-1}\, dt'\\
\frac{4|\Omega|}{TOL^2} \int_{0}^t \left( \sum_{\Delta} \frac{1}{|\Delta|}\int_{\Delta} |\nabla \cdot u^h_\epsilon(t')|^4\, dx \right)\, dt'\leq \|u_0\|^2+\frac{1}{\nu} \int_{0}^t \|f(t')\|^2_{-1}\, dt'\\
\frac{4|\Omega|}{\max |\Delta| TOL^2} \int_{0}^t \left( \sum_{\Delta} \int_{\Delta} |\nabla \cdot u^h_\epsilon(t')|^4\, dx \right)\, dt'\leq \|u_0\|^2+\frac{1}{\nu} \int_{0}^t \|f(t')\|^2_{-1}\, dt'\\
\frac{4|\Omega|}{\max |\Delta| TOL^2} \int_{0}^t  \|\nabla \cdot u^h_\epsilon(t')\|^4_{L^4} \, dt'\leq \|u_0\|^2+\frac{1}{\nu} \int_{0}^t \|f(t')\|^2_{-1}\, dt'.
\end{gathered}
\end{equation*}
The desired result follows.
\end{proof}
\end{proposition}
Let $a(t)=C(\nu)\|\nabla u\|^4+\frac{1}{2}$. We assume $\|\nabla u\|\in L^4(0,T)$, then $a(t) \in L^1(0,T)$. For $0\leq t\leq T$,
\begin{equation}\label{A_t}
    A(t):=\int_{0}^t a(t')\, dt'<\infty.
\end{equation}
Define the error $e =u-u^h_\epsilon$. We give error estimates of the locally adaptive $\epsilon$ penalty method (\ref{penaltyNSE}) to the exact solution of the NSE in the following. 
\begin{theorem}
Let $u$ and $p$ be sufficiently smooth solutions to the NSE, particularly, $\|\nabla u\|\in L^4 (0, T)$, $\|\nabla p\| \in L^2(0, T)$, and $\|p\|\in L^{4/3}(0, T)$. Then there holds
\begin{equation*}
\begin{gathered}
    \sup_{t\in [0,T]}\|e(t)\|^2 + \int_{0}^t \frac{\nu}{4}\|\nabla e\|^2+\sum_{\Delta}\frac{C_1}{LocTol_\Delta}\|\nabla \cdot e\|^4_{L^4(\Delta)}\, dt'
    \leq e^{A(T)} \|e(0)\|^2\\
    +C \inf_{v^h \in X^h} 
    \bigg\{ C\int_{0}^t
 h^2\|(u- v^h)_t\|^2+ \|\nabla  (u-v^h) \|^2 +h^{-2}\|u-v^h \|^2\, dt'\\
    +C \|\nabla (u-v^h)\|^2_{L^4(0,T;L^2)} + \sup_{t\in [0,T]}\|u -v^h\|^2
    +C\int_{0}^t 
     \sum_{\Delta} \frac{\|\nabla \cdot v^h\|^4_{L^4(\Delta)}}{LocTol_\Delta}\, dt'  \bigg\}\\
+Ch^2 \int_{0}^t \| \nabla p\|^2\, dt'+ C \left(\frac{\max|\Delta|}{2}\right)^{1/4} TOL^{1/2}  \left(\int_{0}^t \|p\|^{4/3}\, dt'\right)^{3/4},
\end{gathered}
\end{equation*}
where the constant $C >0$ is dependent on $u_0, \nu, T,f$ and $\Omega$.
\end{theorem}
\begin{remark}
If $X^h$ has a divergence-free subspace with good approximation properties,
\begin{equation*}
\int_{0}^T \sum_{\Delta} \frac{\|\nabla \cdot v^h\|^4_{L^4(\Delta)}}{LocTol_\Delta}\, dt'
\end{equation*}
can vanish by choosing $v^h$ as an approximation of $u$ in that subspace. 
\end{remark}
\begin{proof} 
Subtract equation (\ref{penaltyNSE}) from equation (\ref{weakNSE}), and take $v=v^h$:
\begin{equation}
\begin{gathered}
    (e_t, v^h)+(u\cdot \nabla u,v^h)-b^*(u^h_\epsilon,u^h_\epsilon,v^h)-(p,\nabla \cdot v^h)+\nu(\nabla e,\nabla v^h)\\
    -\sum_{\Delta } LocTol_{\Delta}^{-1}  (|\nabla \cdot u^h_\epsilon |^2 \nabla \cdot u^h_\epsilon \nabla \cdot v^h)=0. \label{erroreqn}
\end{gathered}
\end{equation}
Let $\tilde{u}\in X^h$, define $\eta=u-\tilde{u}$ and $\phi^h=u^h_\epsilon-\tilde{u}$. We write $e=\eta-\phi^h$. For convenience, we define
\begin{equation*}
    a(u,v,w):=\sum_{\Delta} LocTol_{\Delta}^{-1} \int_{\Delta} |\nabla \cdot u|^2 \nabla \cdot v \nabla \cdot w\, dx.
\end{equation*}
Set $v^h=\phi^h$. The error equation (\ref{erroreqn}) becomes:
\begin{equation}\label{error-eqn-bound}
    \begin{gathered}
    \frac{1}{2} \frac{d}{dt}\|\phi^h\|^2+ \nu \|\nabla \phi^h\|^2+
    a(u^h_\epsilon, u^h_\epsilon,\phi^h)=\\
    (\eta_t,\phi^h)+
    \nu(\nabla \eta, \nabla \phi^h)-(p,\nabla \cdot \phi^h)+
    (u\cdot \nabla u, \phi^h)-b^*(u^h_\epsilon,u^h_\epsilon,\phi^h).
    \end{gathered}
\end{equation}
We will first analyze the penalty-related term $a(u^h_\epsilon, u^h_\epsilon, \phi^h)$. Then, we will bound the pressure-related term $(p,\nabla \phi^h)$. Finally, we will address the standard terms of the NSE.

{\bf Penalty term}. Since $\nabla \cdot u = 0$, we have $a(u,u,\phi^h)=0$. Add $a(u,u,\phi^h)$ to the right-hand side of equation (\ref{error-eqn-bound}). Then we subtract $a(\Tilde{u},\Tilde{u},\phi^h) $ from both sides of equation (\ref{error-eqn-bound}). Apply Lemma \ref{strongMonotonicity} equation (\ref{StongMonoton}) to $a(u^h_\epsilon, u^h_\epsilon,\phi^h)-a(\tilde{u},\tilde{u},\phi^h)$:
\begin{equation}
\begin{gathered}
    a(u^h_\epsilon, u^h_\epsilon,\phi^h)-a(\tilde{u},\tilde{u},\phi^h)\\
    =\sum_{\Delta}\frac{1}{LocTol_\Delta} \int_{\Delta}(|\nabla \cdot u^h_\epsilon|^2 \nabla \cdot u^h_\epsilon-|\nabla \cdot \tilde{u}|^2 \nabla \cdot \tilde{u}) \nabla \cdot (u^h_\epsilon-\tilde{u})\, dx\\
    \geq \sum_{\Delta}\frac{1}{LocTol_\Delta} C_1 \int_{\Delta} |\nabla \cdot (u^h_\epsilon-\tilde{u})|^4\,dx\\
    =\sum_{\Delta}\frac{1}{LocTol_\Delta} C_1 \|\nabla \cdot \phi^h\|^4_{L^4(\Delta)}.
    \end{gathered}
\end{equation}
Let $r_\Delta=\max\{ \|\nabla \cdot u\|_{L^4(\Delta)}, \|\nabla \cdot \tilde{u}\|_{L^4(\Delta)}\}=\|\nabla \cdot \tilde{u}\|_{L^4(\Delta)}$. Apply Lemma \ref{strongMonotonicity} equation (\ref{Lipschitz}) to $a(u,u,\phi^h)-a(\tilde{u},\tilde{u},\phi^h)$:
\begin{equation*}
    \begin{gathered}
    a(u,u,\phi^h)-a(\tilde{u},\tilde{u},\phi^h)\\
    =\sum_{\Delta} \frac{1}{LocTol_\Delta} \int_{\Delta} \left(|\nabla \cdot u|^2 \nabla \cdot u -|\nabla \cdot \tilde{u}|^2 \nabla \cdot \tilde{u}\right)\nabla \cdot \phi^h\, dx\\
    \leq \sum_{\Delta}\frac{1}{LocTol_\Delta} C_2 r^2_\Delta \left(\int_{\Delta} |\nabla \cdot (u-\tilde{u})|^4 \, dx)^{1/4}(\int_{\Delta} |\nabla \cdot \phi^h|^4\,dx\right)^{1/4}\\
    =\sum_{\Delta}\frac{1}{LocTol_\Delta} C_2 r^2_\Delta \|\nabla \cdot \eta \|_{L^4(\Delta)} \|\nabla \cdot \phi^h\|_{L^4(\Delta)},
    \end{gathered}
\end{equation*}
Since $\nabla \cdot u =0$, we have $r_\Delta=\|\nabla\cdot (u-\tilde{u})\|_{L^4(\Delta)}=\|\nabla \cdot \eta\|_{L^4(\Delta)}
$.
Thus,
\begin{equation*}
    \begin{split}
    a(u,u,\phi^h)-a(\tilde{u},\tilde{u},\phi^h)\leq
    \sum_{\Delta}\frac{C_2}{LocTol_\Delta} \|\nabla \cdot \eta \|^3_{L^4(\Delta)} \|\nabla \cdot \phi^h\|_{L^4(\Delta)}
    \end{split}.
\end{equation*}
Apply Lemma (\ref{holderandyoung}) H\"older's and Young's inequality with $p=4$ and $q=\frac{4}{3}$ to the following term:
\begin{equation*}
\begin{gathered}
    \sum_{\Delta}\frac{C_2}{LocTol_\Delta} \|\nabla \cdot \eta \|^3_{L^4(\Delta)} \|\nabla \cdot \phi^h\|_{L^4(\Delta)}\\
    =\sum_{\Delta} \left( \frac{C_1^{1/4}}{LocTol^{1/4}_\Delta} \|\nabla \cdot \phi^h\|_{L^4(\Delta)}\right)\left(\frac{C_2}{C_1^{1/4}LocTol^{3/4}_\Delta} \|\nabla \cdot \eta \|^3_{L^4(\Delta)} \right) \\
    \leq \left(\sum_{\Delta}  \frac{C_1}{LocTol_\Delta} \|\nabla \cdot \phi^h\|^4_{L^4(\Delta)}\right)^{1/4} \left(\sum_{\Delta}\frac{C_2^{4/3}}{C_1^{1/3}LocTol_\Delta} \|\nabla \cdot \eta \|^{4}_{L^4(\Delta)} \right)^{3/4}\\
    \leq \frac{1}{2}\sum_{\Delta} \frac{C_1}{LocTol_\Delta} \|\nabla \cdot \phi^h\|^4_{L^p(\Delta)} +\frac{2^{-1/3}}{4/3}\sum_{\Delta}\frac{C_2^{4/3}}{C_1^{1/3}LocTol_\Delta} \|\nabla \cdot \eta \|^{4}_{L^4(\Delta)}.
    \end{gathered}
\end{equation*}
We have
\begin{equation}
\begin{gathered}
a(u,u,\phi^h)-a(\tilde{u},\tilde{u},\phi^h)\leq \frac{1}{2}\sum_{\Delta} \frac{C_1}{LocTol_\Delta} \|\nabla \cdot \phi^h\|^4_{L^p(\Delta)} \\
+\frac{2^{-1/3}}{4/3}\sum_{\Delta}\frac{C_2^{4/3}}{C_1^{1/3}LocTol_\Delta} \|\nabla \cdot \eta \|^{4}_{L^4(\Delta)}.
\end{gathered}
\end{equation}
We will hide the term $\frac{1}{2}\sum_{\Delta} \frac{C_1}{LocTol_\Delta} \|\nabla \cdot \phi^h\|^4_{L^p(\Delta)}$ on the left.

{\bf Pressure related term.}
The proof to bound $(p,\nabla \phi^h)$ is inspired by Falk \cite{falk1976finite}. We have
\begin{equation*}
\begin{split}
    (p,\nabla \cdot \phi^h)
    =(p,\nabla \cdot u^h_\epsilon)+(p,\nabla \cdot (u-\tilde{u}))=(p,\nabla \cdot u^h_\epsilon)+(p,\nabla \cdot \eta).
\end{split}
\end{equation*}
We bound  $(p,\nabla \cdot u^h_\epsilon)$ and $(p,\nabla \cdot \eta)$ separately. To bound $(p,\nabla \cdot \eta)$, we use integration by parts and Lemma \ref{holderandyoung}:
\begin{equation*}
\begin{split}
(p,\nabla \cdot \eta)
    =-(\nabla p, \eta)+ \int_{\partial \Omega} p (\eta \cdot n)\,ds
    \leq \|\nabla p\|\|\eta\|
    \leq \frac{h^2}{2}\|\nabla p\|^2+\frac{1}{2h^2}\|\eta\|^2.
    \end{split}
\end{equation*}
For the term $(p,\nabla \cdot u^h_\epsilon)$, we localize $\nabla \cdot u^h_\epsilon$ and use Lemma \ref{holderandyoung} and Cauchy-Schwarz inequality inequality. The proof follows Xie \cite{xie2022adaptive} p.8. So it is omitted here. There results
\begin{equation*}
\begin{split}
        (p,\nabla \cdot u^h_\epsilon)
    \leq \left(\sum_{\Delta} \frac{\|\nabla \cdot u^h_\epsilon\|^4_{L^4(\Delta)}}{LocTol_\Delta}\right)^{1/4} \left(\frac{\max|\Delta|}{2}\right)^{1/4} TOL^{1/2} \|p\|.
    \end{split}
\end{equation*}
Thus, we have
\begin{equation}
\begin{gathered}
    (p,\nabla \cdot \phi^h)\leq \left(\sum_{\Delta} \frac{\|\nabla \cdot u^h_\epsilon\|^4_{L^4(\Delta)}}{LocTol_\Delta}\right)^{1/4} \left(\frac{\max|\Delta|}{2}\right)^{1/4} TOL^{1/2} \|p\|\\
    +\frac{h^2}{2}\|\nabla p\|^2+\frac{1}{2h^2}\|\eta\|^2.
\end{gathered}
\end{equation}

{\bf Standard NSE terms.} These are bounded by adapting techniques in e.g. John \cite{john2016finite} and Layton \cite{layton2008introduction}.
Using Lemma \ref{holderandyoung} H\"older's and Young's inequality, we have
\begin{equation}
\begin{split}
(\eta_t, \phi^h)\leq \|\eta_t\|\| \phi^h\|\leq \frac{1}{2}\|\phi^h\|^2 + \frac{1}{2}\|\eta_t\|^2. 
\end{split}
\end{equation}
\begin{equation}
\begin{split}
\nu(\nabla \eta, \nabla \phi^h)\leq \|\nabla \eta\|\|\nabla \phi^h\|\leq \frac{\nu}{8}\|\nabla \phi^h\|^2 + C(\nu) \|\nabla \eta\|^2.
\end{split}
\end{equation}
Last, we work on the trilinear terms $(u\cdot \nabla u, \phi^h)- b^*(u^h_\epsilon, u^h_\epsilon, \phi^h)$. Since $\nabla \cdot u=0$, $(u\cdot \nabla u,\phi^h) = b^*(u, u,\phi^h)$. Thus, it is equivalent to bound $b^*(u,u,\phi^h)-b^*(u^h_\epsilon,u^h_\epsilon,\phi^h)$. We use add and subtract techniques and $e= \eta -\phi^h$. 
\begin{equation*}
    \begin{gathered}
    b^*(u,u,\phi^h)- b^*(u^h_\epsilon,u^h_\epsilon,\phi^h)\\
        =b^*(u,u,\phi^h)-b^*(u^h_\epsilon,u,\phi^h)+b^*(u^h_\epsilon,u,\phi^h)-b^*(u^h_\epsilon,u^h_\epsilon,\phi^h)\\
    =b^*(\eta,u,\phi^h)-b^*(\phi^h,u,\phi^h)+b^*(u^h_\epsilon,\eta,\phi^h)-b^*(u^h_\epsilon,\phi^h,\phi^h).
    \end{gathered}
\end{equation*}
Since $b^*(u^h_\epsilon,\phi^h,\phi^h)=0$, we have
\begin{equation*}
    \begin{split}
    b^*(u,u,\phi^h)- b^*(u^h_\epsilon,u^h_\epsilon,\phi^h)
        =b^*(\eta,u,\phi^h)-b^*(\phi^h,u,\phi^h)+b^*(u^h_\epsilon,\eta,\phi^h).
    \end{split}
\end{equation*}
We use Lemma \ref{upperbounds} and Lemma \ref{holderandyoung} to bound the trilinear terms. 
\begin{equation}
\begin{split}
    b^*(\eta,u,\phi^h)
    \leq C \|\nabla \eta\| \|\nabla u\|\|\nabla \phi^h\|\leq \frac{\nu}{8}\| \nabla \phi^h\|^2+ C(\nu) \|\nabla u\|^2 \|\nabla \eta\|^2,\\
\end{split}
\end{equation}
\begin{equation}
\begin{gathered}
    b^*(\phi^h, u,\phi^h)\leq C\|\phi^h\|^{1/2}\|\nabla \phi^h\|^{3/2}\|\nabla u\|\\
    \leq \frac{\nu}{8}\|\nabla \phi^h\|^2+C(\nu)\|\nabla u\|^4\|\phi^h\|^2,\\
        \end{gathered}
\end{equation}
\begin{equation*}
\begin{gathered}
    b^*(u^h_\epsilon,\eta,\phi^h)\leq  C\|u^h_\epsilon\|^{1/2}\|\nabla u^h_\epsilon\|^{1/2}\|\nabla \eta\|\|\nabla \phi^h\|\\
    \leq \frac{\nu}{8} \|\nabla \phi^h\|^2 + C(\nu)\|u^h_\epsilon\|\|\nabla u^h_\epsilon\|\|\nabla \eta\|^2.
    \end{gathered}
\end{equation*}
By Theorem \ref{ppstability}, the following holds:
\begin{equation*}
\sup_{0\leq t\leq T} \|u^h_\epsilon\|^2+\|\nabla u^h_\epsilon\|^2_{L^2(0,T;L^2)}\leq \|u_0\|^2+\frac{1}{\nu}\int_{0}^T \|f(t')\|^2_{-1}\, dt'.
\end{equation*}
We have
\begin{equation}
    b^*(u^h_\epsilon,\eta,\phi^h)\leq  \frac{\nu}{8} \|\nabla \phi^h\|^2 + C\|\nabla u^h_\epsilon\|\|\nabla \eta\|^2.
\end{equation}
Thus, we have
\begin{equation}\label{before-integrating-factor}
    \begin{gathered}
    \frac{1}{2} \frac{d}{dt}\|\phi^h\|^2+ \frac{\nu}{2} \|\nabla \phi^h\|^2+
    \frac{1}{2}\sum_{\Delta}\frac{C_1}{LocTol_\Delta}  \|\nabla \cdot \phi^h\|^4_{L^4(\Delta)}
    \leq \frac{1}{2}\|\eta_t\|^2+\\
    C(\nu)\|\nabla \eta\|^2+\left(\sum_{\Delta} \frac{\|\nabla \cdot u^h_\epsilon\|^4_{L^4(\Delta)}}{LocTol_\Delta}\right)^{1/4} \left(\frac{\max|\Delta|}{2}\right)^{1/4} TOL^{1/2} \|p\|\\
    +\frac{h^2}{2}\|\nabla p\|^2+\frac{\|\eta\|^2}{2h^2}
    +C\|\nabla u^h_\epsilon\| \|\nabla \eta\|^2
     +C(\nu)\|\nabla u\|^2 \|\nabla \eta\|^2\\
     +\sum_{\Delta}\frac{C(C_1,C_2)}{LocTol_\Delta} \|\nabla \cdot \eta \|^{4}_{L^4(\Delta)}
    +\left(C(\nu) \|\nabla u\|^4 +\frac{1}{2}\right)\|\phi^h\|^2.
\end{gathered} 
\end{equation}
Multiply both sides of equation (\ref{before-integrating-factor}) by the integration factor $e^{-A(t)}$, where $A(t)$ is defined in equation (\ref{A_t}). This gives
\begin{equation*}
\begin{gathered}
    \frac{d}{dt}\left( e^{-A(t)}\|\phi^h\|^2\right)+e^{-A(t)}\left(\frac{\nu}{4}\|\nabla \phi^h\|^2+\sum_{\Delta}\frac{C_1}{LocTol_\Delta} \|\nabla \cdot \phi^h\|^4_{L^4(\Delta)}\right) \\
    \leq e^{-A(t)}h^2\|\nabla p\|^2+2e^{-A(t)}\sum_{\Delta} \frac{C(C_1,C_2)}{LocTol_\Delta}\|\nabla \cdot \eta\|^4_{L^4(\Delta)}
    \\
    + e^{-A(t)} C
    \left( \|\eta_t\|^2+ \left(1+\|\nabla u\|^2+ \|\nabla u^h_\epsilon\|\right)\|\nabla \eta\|^2+\frac{\|\eta\|^2 }{h^2}\right)\\
    +2e^{-A(t)}\left(\sum_{\Delta} \frac{\|\nabla \cdot u^h_\epsilon\|^4_{L^4(\Delta)}}{LocTol_\Delta}\right)^{1/4} \left(\frac{\max|\Delta|}{2}\right)^{1/4} TOL^{1/2} \|p\|. 
    \end{gathered}
\end{equation*}
Integrate from $0$ to $t$, and multiply $e^{A(t)}$ to both sides. We have
\begin{equation*}
\begin{gathered}
    \|\phi^h(t)\|^2 + \int_{0}^t e^{A(t)-A(t')}\left(\frac{\nu}{4}\|\nabla \phi^h\|^2+\sum_{\Delta}\frac{C_1}{LocTol_\Delta} \|\nabla \cdot \phi^h\|^4_{L^4(\Delta)}\right)\, dt'\\
    \leq e^{A(t)} \|\phi^h(0)\|^2+\int_{0}^t e^{A(t)-A(t')} h^2\|\nabla p\|^2\, dt'\\
+\int_{0}^t 2 e^{A(t)-A(t')}
     \sum_{\Delta} \frac{C(C_1,C_2)}{LocTol_\Delta}\|\nabla \cdot \eta\|^4_{L^4(\Delta)}\, dt'\\
+C\int_{0}^t e^{A(t)-A(t')}
    \left( \|\eta_t\|^2+\left(1+\|\nabla u\|^2+ \|\nabla u^h_\epsilon\|\right)\|\nabla \eta\|^2 +h^{-2}\|\eta\|^2 \right)\, dt'\\
+\int_{0}^t 2e^{A(t)-A(t')}\left(\sum_{\Delta} \frac{\|\nabla \cdot u^h_\epsilon\|^4_{L^4(\Delta)}}{LocTol_\Delta}\right)^{1/4} \left(\frac{\max|\Delta|}{2}\right)^{1/4} TOL^{1/2} \|p\|  \, dt'.
\end{gathered}
\end{equation*}
Since $t'\leq t$, $A(t)-A(t')\geq 0$,  $e^{A(t)-A(t')}\geq 1$. $A(t)-A(t')\leq A(T)$ for all $ t'\leq t\leq T$. We have
\begin{equation}\label{error-phi}
\begin{gathered}
    \|\phi^h(t)\|^2 + \int_{0}^t \left(\frac{\nu}{4}\|\nabla \phi^h\|^2+\sum_{\Delta}\frac{C_1}{LocTol_\Delta} \|\nabla \cdot \phi^h\|^4_{L^4(\Delta)}\right)\, dt'\\
    \leq e^{A(T)} \|\phi^h(0)\|^2+e^{A(T)}\int_{0}^t  h^2\|\nabla p\|^2\, dt'
    \\
+e^{A(T)}\int_{0}^t 2 
     \sum_{\Delta} \frac{C(C_1,C_2)}{LocTol_\Delta}\|\nabla \cdot \eta\|^4_{L^4(\Delta)}\, dt'\\
     +Ce^{A(T)}\int_{0}^t 
     \|\eta_t\|^2+\left(1+\|\nabla u\|^2+ \|\nabla u^h_\epsilon\|\right)\|\nabla \eta\|^2+\frac{\|\eta\|^2}{h^2} \, dt'\\
    +2e^{A(T)}\int_{0}^t \left(\sum_{\Delta} \frac{\|\nabla \cdot u^h_\epsilon\|^4_{L^4(\Delta)}}{LocTol_\Delta}\right)^{1/4} \left(\frac{\max|\Delta|}{2}\right)^{1/4} TOL^{1/2} \|p\| \, dt'.
\end{gathered}
\end{equation}
We can simplify equation (\ref{error-phi}) by the following estimates.
\begin{equation*}
    \begin{gathered}
        \int_{0}^t \|\nabla u\|^2 \|\nabla \eta\|^2\, dt' \leq \left(\int_{0}^t  \|\nabla u\|^4 \, dt' \right)^{1/2} \left(\int_{0}^t \|\nabla \eta\|^4\, dt'\right)^{1/2}\\
        \leq \|\nabla u\|^2_{L^4(0,T;L^2)}\|\nabla \eta\|^2_{L^4(0,T;L^2)}. \\
        \end{gathered}
        \end{equation*}
\begin{equation*}
\begin{gathered}
        \int_{0}^t \|\nabla u^h_\epsilon\| \|\nabla \eta\|^2\, dt'\leq \left(\int_{0}^t  \|\nabla u^h_\epsilon\|^2 \, dt' \right)^{1/2} \left(\int_{0}^t \|\nabla \eta\|^4\, dt'\right)^{1/2}\\
        \leq \|\nabla u^h_\epsilon\|_{L^2(0,T;L^2)}\|\nabla \eta\|^2_{L^4(0,T;L^2)}.
    \end{gathered}
\end{equation*}
By Theorem \ref{stability-bound}, $\|\nabla u^h_\epsilon\|^2_{L^2(0,T;L^2)}$ is bounded by data. Apply Lemma \ref{holderandyoung} and Theorem \ref{stability-bound} to the following term:
\begin{equation*}
\begin{gathered}
    \int_{0}^t\left(\sum_{\Delta} \frac{\|\nabla \cdot u^h_\epsilon\|^4_{L^4(\Delta)}}{LocTol_\Delta}\right)^{1/4} \left(\frac{\max|\Delta|}{2}\right)^{1/4} TOL^{1/2} \|p\|  \, dt'\\
    =\left(\frac{\max|\Delta|}{2}\right)^{1/4} TOL^{1/2} \int_{0}^t\left(\sum_{\Delta} \frac{\|\nabla \cdot u^h_\epsilon\|^4_{L^4(\Delta)}}{LocTol_\Delta}\right)^{1/4}  \|p\|  \, dt'\\
    \leq \left(\frac{\max|\Delta|}{2}\right)^{1/4} TOL^{1/2} \left(\int_{0}^t \sum_{\Delta} \frac{\|\nabla \cdot u^h_\epsilon\|^4_{L^4(\Delta)}}{LocTol_\Delta}\, dt'\right)^{1/4} \left(\int_{0}^t \|p\|^{4/3}\, dt'\right)^{3/4} \\
    \leq C\left(\frac{\max|\Delta|}{2}\right)^{1/4} TOL^{1/2} \left(\int_{0}^t \|p\|^{4/3}\, dt'\right)^{3/4}.
\end{gathered}
\end{equation*}
Thus,
\begin{equation*}
\begin{gathered}
    \|\phi^h(t)\|^2 + \int_{0}^t \left(\frac{\nu}{4}\|\nabla \phi^h\|^2+\sum_{\Delta}\frac{C_1}{LocTol_\Delta} \|\nabla \cdot \phi^h\|^4_{L^4(\Delta)}\right)\, dt'\\
    \leq e^{A(T)} \|\phi^h(0)\|^2+ e^{A(T)}h^2 \int_{0}^t  \|\nabla p\|^2\, dt' \\
    +2 e^{A(T)}\int_{0}^t 
     \sum_{\Delta} \frac{C(C_1,C_2)}{LocTol_\Delta}\|\nabla \cdot \eta\|^4_{L^4(\Delta)}\, dt'\\
    +C\int_{0}^t
    \left( \|\eta_t\|^2+\|\nabla \eta\|^2+h^{-2}\|\eta\|^2 \right)\, dt' + C\|\nabla \eta\|^2_{L^4(0,T;L^2)}\\
    +Ce^{A(T)} \left(\frac{\max|\Delta|}{2}\right)^{1/4} TOL^{1/2} \left(\int_{0}^t \|p\|^{4/3}\, dt'\right)^{3/4}.
\end{gathered}
\end{equation*}
We use the triangle inequality $\|e\|\leq \|\eta\|+ \|\phi^h\|$ to have the final result.
\end{proof}
\begin{proposition}
If $X^h$ has a divergence-free subspace with good approximations,
\begin{equation*}
\begin{gathered}
    \sup_{t\in [o,T]}\|u(t)-u^h_\epsilon(t)\|^2 +\int_{0}^T\frac{\nu}{4}\|\nabla u^h_\epsilon(t)\|^2\, dt'\\
    +\int_{0}^T \sum_{\Delta}\frac{C_1 \|\nabla \cdot (u(t)-u^h_\epsilon(t))\|^4_{L^4(\Delta)}}{LocTol_\Delta} \, dt'\\
    \leq e^{A(T)} \|u(0)-u^h_\epsilon(0)\|^2 + C h^{2m+2} \sup_{t\in [0,T]}|u|_{m+1}\\
    +C h^{2m} \int_{0}^T
   |u_t|^2_{m+1}+ |u|^2_{m+1} \, dt'+ C h^{2m} \left(\int_{0}^T |u|^4_{m+1}\, dt' \right)^2 \\
+Ch^2 \int_{0}^T  \|\nabla p\|^2\, dt'+ C \left(\frac{\max|\Delta|}{2}\right)^{1/4} TOL^{1/2}  \left(\int_{0}^T \|p\|^{4/3}\, dt'\right)^{3/4}.
\end{gathered}
\end{equation*} 
\end{proposition} 
\begin{remark} If we pick $\tilde{u}$ the $L^2$ projection of $u$, then by Lemma \ref{h-inverse-projection}, the term $\|\eta_t\|^2$ can be replaced by $h^2\|\eta_t\|^2$. Indeed, in this case, $\|\eta_t\|^2_{-1}\leq h^2 \|\eta_t\|^2$, and we have 
\begin{equation}
\begin{split}
(\eta_t, \phi^h)\leq \|\eta_t\|_{-1}\|\nabla \phi^h\|\leq \frac{\nu}{4}\|\nabla \phi^h\|^2 + C(\nu)\|\eta_t\|^2_{-1}. 
\end{split}
\end{equation}
\end{remark}
\section{Numerical tests}\label{sec: tests} We first use a modified Green-Taylor vortex to calculate convergence rates of the divergence $\|\nabla \cdot u\|$ in terms of the global tolerance $TOL$. We also monitor the corresponding spatial average of the penalty parameter $\epsilon_{ave}$ and velocity errors with different $TOL$. We employ the algorithm on a complex 2d flow in the second test. We compare how well the locally adaptive penalty method (Algorithm \ref{alg: elementwise}) is preserving the divergence-free condition. In addition, we compare it with the constant penalty $\epsilon = \Delta t$ and the non-penalized coupled NSE. We also combine the locally adaptive $\epsilon$ penalty method with a first-order adaptive time stepping (Algorithm \ref{alg: doubly-first-order}) and test it with a flow problem with sharp transition regions. We compare it with the locally adaptive $\epsilon$ penalty method with a constant timestep.     

We use the backward Euler time discretization for the momentum equation in the tests. We use the Taylor-Hood $\left( P2- P1\right)$ finite element pair for approximating the velocity and pressure fields. We generate the unstructured meshes by GMSH \cite{geuzaine2009gmsh} with a target mesh size parameter $lc$. 
\subsection{Test 1: the modified Green-Taylor vortex for convergence \cite{fang2023penalty}}
\label{subsec:2}
The Green-Taylor vortex problem is often employed to calculate convergence rates due to its known analytical solution \cite{john2016finite}. In our adaptation of the Green-Taylor experiment within the domain $\Omega = (0,1)^2$, we have adjusted the velocity to ensure it does not diminish. The analytical solution for the modified Green-Taylor vortex in this context is as follows:
\begin{align*}
\begin{gathered}
    u(x,y,t) = -cos(x) sin(y) sin(t),\\
    v(x,y,t) = sin(x) cos(y) sin(t),\\
    p(x,y,t) = \frac{1}{4}\left( cos(2x)+ cos(2y)\right) sin^2(t).
\end{gathered}
\end{align*}
We take $\nu = 1, \mathcal{R}e = 1/\nu,$ mesh size $h = lc = 1/27$, $\Delta t = h^2$, and the final time $T=1$. The range of the penalty parameter is from  $EpsMin=1\times 10^{-6}$  to  $EpsMax = 1\times 10^{-1}$. We impose the exact solution on the boundary. We calculate some flow statistics at time $T=1$ as $TOL$ changes, including velocity errors, the divergence of velocity, and the spatial average of the penalty parameter $\epsilon_{ave}$. Denote $\|\nabla \cdot u\| = C\ {TOL}^\alpha$, where $C$ is a constant. Solve $\alpha$ via 
\begin{equation*}
\alpha = \frac{ln(\|\nabla \cdot u\|_{TOL_1}/\|\nabla \cdot u\|_{TOL_2})}{ln(TOL_1/TOL_2)}.
\end{equation*}  
\begin{table}[H]
\caption{$\|\nabla \cdot u\| = \mathcal{O}(TOL)$}
\label{tab:1}
\centering
\begin{tabular}{c |c| c | c | c | c }
$TOL$ & $\|\nabla \cdot u^h_\epsilon\|$ & $\alpha$ & $\epsilon_{ave}$& $\max_{t_n}\|u-u^h_\epsilon\|$ & $\int_{0}^T\|\nabla u-\nabla u^h_\epsilon\|$  \\
\hline
$ 10^{-1} $ & $8.7\times 10^{-3}$ &-- & $1.0\times 10^{-1}$&  $2.9\times 10^{-3}$ & $7.0 \times 10^{-3}$ \\ 
$10^{-2} $ & $5.7\times 10^{-3}$ & 0.18& $8.1\times 10^{-2}$ & $1.8 \times 10^{-3}$ & $5.4\times 10^{-3}$ \\
$10^{-3} $ & $7.0\times 10^{-4}$&0.91& $1.8\times 10^{-2}$  & $2.0 \times 10^{-4}$ &$8.6 \times 10^{-4}$\\
$10^{-4} $ & $7.1\times 10^{-5}$&0.99& $3.4\times 10^{-4}$& $2.0 \times 10^{-5}$ &$6.7 \times 10^{-4}$ \\
$10^{-5} $& $6.2\times 10^{-5}$&0.06& $1.0\times 10^{-6}$ & $9.3 \times 10^{-6}$ & $8.6 \times 10^{-4}$\\ 
\end{tabular}
\end{table}
As stated in Table \ref{tab:1}, we observe first-order convergence for the divergence $\|\nabla \cdot u^h_\epsilon\| = \mathcal{O}(TOL)$. Proposition \ref{l4-div-bound} suggests $\int_{0}^t \|\nabla \cdot u^h_\epsilon(t')\|^4_{L^4}\, dt' =\mathcal{O}(TOL^2)$. When $TOL=10^{-5}$, the $\epsilon_{ave}$ reaches the lower bound of user-defined $EpsMin$. Velocity errors get smaller as we set smaller $TOL$ until they saturate at $\mathcal{O}(TOL)$.
\subsection{Test 2: flows between offset cylinders taken from \cite{layton2020doubly}, \cite{xie2022adaptive}}
Consider the flow between offset circles to test how well the incompressibility condition is satisfied using the locally adaptive penalty method. The domain is a disk with a smaller offset obstacle inside. $$\Omega = \{(x,y): x^2 + y^2 \leq r_1^2 \cap (x-c_1)^2 + (y-c_2)^2 \geq r_2^2 \}. $$ where $r_1 =1, r_2 = 0.1, c = (c_1,c_2)=(\frac{1}{2},0)$. We generate the unstructured meshes with GMSH \cite{geuzaine2009gmsh}, with a target mesh size $lc=0.01$. The area of the triangle elements ranges from $2.2\times 10^{-5}$ to $6.0 \times 10^{-5}$. The flow is driven by a counterclockwise force $f(x,y,t) = (4x\min \left(t,1\right)(1-x^2 -y^2) , -4y\min \left(t,1\right) (1-x^2-y^2))$. Flow is at rest at $t=0$. Impose the no-slip boundary conditions on both circles. We set $\Delta t =0.02$, $\nu = 0.01$, $L = 1$, $U = 1$ and $\text{Re} = \frac{UL}{\nu}$. The final time $T=16$. 
We apply the second order linear extrapolation of $u$ to the nonlinear term, $b^*(u^*, u^{h}_{\epsilon, n+1}, v^h)$. For constant timestep it is
\begin{equation*}
u^*=2u^{h}_{\epsilon,n} -u^{h}_{\epsilon,n-1}.
\end{equation*}
We add a time filter to increase the accuracy from first-order to second-order  \cite{kean2023doubly,guzel2018time}:
\begin{equation*}
\begin{gathered}
    (\frac{u^{h,1}_{\epsilon,n+1}-u^h_{\epsilon,n}}{\Delta t},v^h)+b^*(u^*,u^{h,1}_{\epsilon,n+1},v^h)+\nu (\nabla u^{h,1}_{\epsilon,n+1},\nabla v^h)\\
    +\sum_{\Delta} \epsilon_{\Delta}^{-1}\int_{\Delta}  \nabla \cdot u^{h,1}_{\epsilon,n+1} \nabla \cdot v^h\, dx
    =(f^{n+1},v^h),
    \end{gathered}
\end{equation*}
\begin{equation*}
    u^{h}_{\epsilon,n+1} = u^{h,1}_{\epsilon,n+1} -\frac{1}{3} (u^{h,1}_{\epsilon,n+1}-2u^{h}_{\epsilon,n}+ u^{h}_{\epsilon,n-1}).
\end{equation*}
\begin{figure}
    \centering
\includegraphics[width=0.85\textwidth]{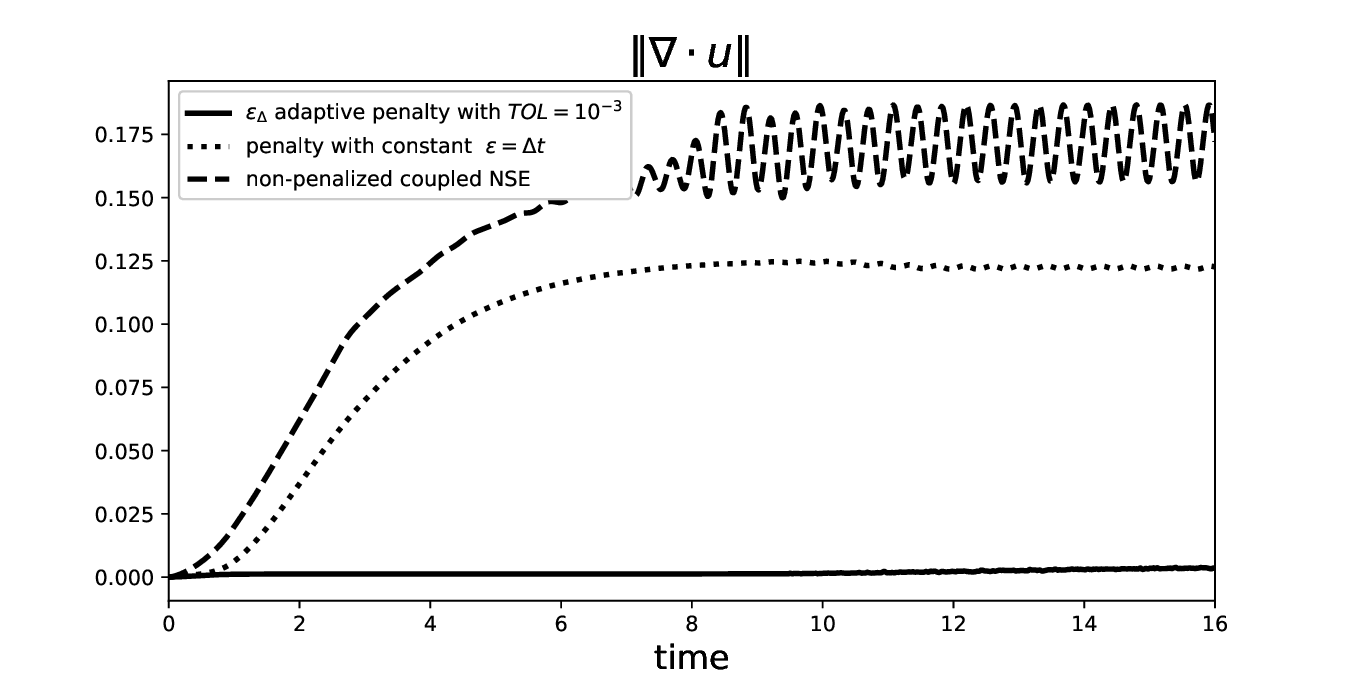}
    \caption{Comparison of the divergence of velocity $\|\nabla \cdot u\|$ for  adaptive $\epsilon$ penalty, constant penalty $\epsilon=\Delta t$, and non-penalized coupled NSE. The adaptive penalty method is far better.}
    \label{fig: div}
\end{figure}
\begin{figure}
    \centering    
    \includegraphics[width=0.85\textwidth]{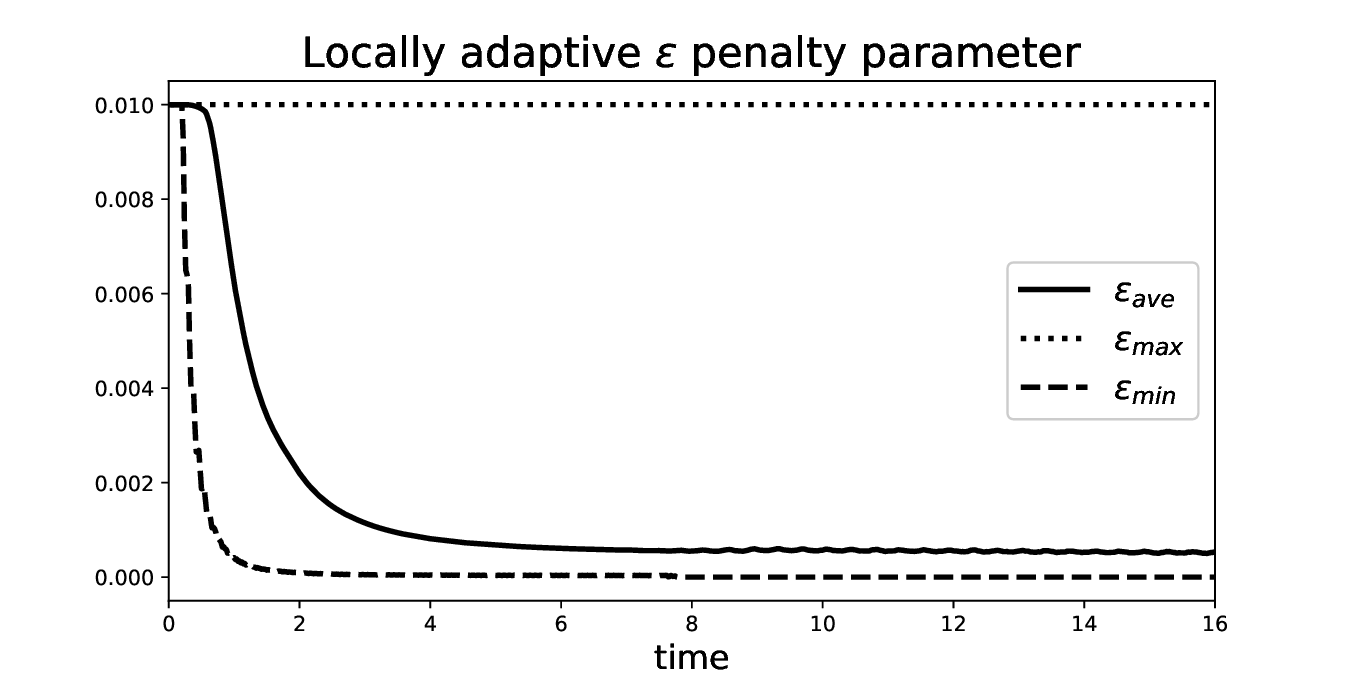}
    \caption{The $\epsilon_{min}$, $\epsilon_{max}$ and $\epsilon_{ave}$ among all the triangle meshes evolve for the adaptive $\epsilon$ penalty. After reaching a steady state, $\epsilon_{min}=EpsMin=10^{-10}$, $\epsilon_{max}=EpsMax=10^{-2}$ and $\epsilon_{ave} \approx 5 \times 10^{-4}$.}
    \label{fig: epsilon}
\end{figure}

Figure \ref{fig: div} illustrates the comparison of the divergence of the velocity, $\|\nabla \cdot u \|$, using three different methods: locally adaptive $\epsilon$ penalty, the constant penalty with $\epsilon=\Delta t=0.02$, and the non-penalized coupled NSE over an extended period. Set the global tolerance $TOL=10^{-3}$ for the locally adaptive penalty method. We observe $\|\nabla \cdot u^h_{\epsilon}\| \leq  4.0 \times10^{-3}$. This algorithm effectively controls the velocity divergence, ensuring small $\|\nabla \cdot u\| $. For the constant penalty case with $\epsilon=\Delta t$, the divergence, $\|\nabla \cdot u^h\|$, is approximately $0.12$, much bigger. The non-penalized coupled NSE does not satisfy the incompressibility condition.

Figure \ref{fig: epsilon} tracks the minimum $\epsilon_{min}$, the maximum $\epsilon_{max}$ and the average of the penalty parameter $\epsilon_{ave}$ over time. Different triangle elements require different $\epsilon$ to ensure $\|\nabla \cdot u\|_{\Delta}\approx LocTol$. It ranges from $EpsMin=10^{-10}$ to $EpsMax=10^{-2}$. The average $\epsilon$ reach a steady state where $\epsilon_{ave} \approx 5\times 10^{-4}$. The locally adaptive penalty algorithm seeks the optimal $\epsilon$ to satisfy the $LocTol_{\Delta}$ at each time step, thus the $TOL$.
\subsection{Locally adaptive penalty with time adaptivity}\label{sec: doubly-adaptive} Adaptive time stepping is necessary when the fluid flow changes rapidly over time. We extend the locally adaptive penalty method by combining it with adaptive time stepping. Kean, Xie, and Xu \cite{kean2023doubly} gave a comprehensive study on the adaptive penalty method (with $\epsilon$ constant in space) with adaptive time steps and variable order of the methods. We evaluate the locally adaptive $\epsilon$ algorithm with adaptive time stepping. In particular, we will combine the locally adaptive algorithm with adaptive time stepping for the first-order method in Algorithm \ref{alg: doubly-first-order}. We compare it with the locally adaptive algorithm.

The doubly adaptive methods of Algorithm \ref{alg: doubly-first-order} build on Kean, Xie, and Xu \cite{kean2023doubly} and Guzel and Layton \cite{guzel2018time} and Layton and McLaughlin \cite{layton2020doubly}. The main idea behind these methods is summarized below. The time discretization from \cite{guzel2018time}: for $y'= f(t,y)$, select $\tau = k_{n+1}/k_n$, $\alpha = \tau (1+\tau)/(1+2\tau)$, where $k_n$ is the time step at $t_n$. Then
\begin{equation*}
\begin{gathered}
\frac{y^1_{n+1}-y_n}{k_{n+1}} = f(t_{n+1}, y_{n+1}),\\
y_{n+1} = y^1_{n+1} -\frac{\alpha}{2} \left( \frac{2 k_n}{k_n+ k_{n+1}} y_{n+1} - 2 y_n + \frac{2 k_{n+1}}{k_n+k_{n+1}} y_{n-1}\right),\\
EST = |y_{n+1}-y^1_{n+1}|.
\end{gathered}
\end{equation*}
Denote $D_2$ the difference at $t_{n+1}$ is given by:
\begin{equation*}
D_2(n+1)= \frac{2k_n}{k_n + k_{n+1}} u^1_{n+1} -2u_n + \frac{2k_{n+1}}{k_n +k_{n+1}} u_{n-1}.
\end{equation*}
The estimate of the local truncation error in the first-order estimation is
\begin{equation*}
tEST_1 = \|u_{n+1}-u^1_{n+1}\| = \frac{\alpha}{2}\|D_2(n+1)\|. 
\end{equation*}
The first-order prediction of the time stepping is
\begin{equation*}
k_{new} = k_{old} \left(\frac{tTOL}{tEST_1}\right)^{1/2}.
\end{equation*}
The basic concept of adaptive time stepping involves increasing the timestep when $tEST_1 < tTOLmin$ and decreasing it when $tEST_1 > tTOL$, repeating the calculation at the current timestep. For detailed time stepping procedures, refer to Algorithm \ref{alg: doubly-first-order}.
\begin{algorithm}
\caption{Locally adaptive $\epsilon$ penalty with variable
time stepping}\label{alg: doubly-first-order}
\begin{algorithmic}
\State Given $tTOL$, $tTOLmin$, $MaxRetry$, $TOL$, initial timestep $\Delta t$, $EpsMin$, and $EpsMax$. Set initial $\epsilon_{\Delta} =1$. Compute the local tolerance on each triangle 
$LocTol_{\Delta}=\frac{1}{2}\frac{TOL^2}{|\Omega|}|\Delta|$.
\While{$t < T$}:
\State Update $t_{n+1}=t_n+ k_{n+1}$.
\State Compute $\tau = \frac{k_{n+1}}{k_n}$, $\alpha_1 =\frac{\tau(1+\tau)}{1+2\tau}$.
\State Compute $u^* = (1+ \tau) u_n - \tau u_{n-1}$.
\State Given $u^h_{\epsilon,n}$, solve for $u^h_{\epsilon,n+1}$ using penalty method: find $u^h_{\epsilon,n+1} \in X^h$ such that
\begin{equation*}
\begin{gathered}
(\frac{u^h_{\epsilon,n+1}-u^h_{\epsilon,n}}{k_{n+1}},v^h)+b^*(u^*,u^h_{\epsilon,n+1},v^h)+\nu (\nabla u^h_{\epsilon,n+1},\nabla v^h)\\
    +\sum_{\Delta} \epsilon_{\Delta}^{-1}\int_{\Delta}  \nabla \cdot u^h_{\epsilon,n+1} \nabla \cdot v^h\, dx=(f^{n+1},v^h).
    \end{gathered}
\end{equation*}
\State Compute $D_2$ and $tEST_1$:
\begin{equation*}
\begin{gathered}
D_2(n+1) = \frac{2k_n}{k_n + k_{n+1}} u^1_{n+1} -2u_n + \frac{2k_{n+1}}{k_n +k_{n+1}} u_{n-1},\\
tEST_1(n+1) = \frac{\alpha_1}{2}\|D_2(n+1)\|.
\end{gathered}
\end{equation*}
\State {\bf Update the timestep $k$ using the decision tree:}
\If{$tEST_1< tTOLmin$}
\begin{equation*}
\begin{gathered}
k_{n+2} = \max\left\{\min\left\{0.9 k_{n+1}(tTOL/tEST_1)^{1/2}, 2 k_{n+1}\right\}, 0.5 k_{n+1}\right\},\\
k_{n+2} = \min\left\{ k_{n+2}, \Delta t_{\max}\right\}.
\end{gathered}
\end{equation*}
\ElsIf{$tEST_1>tTOL$}
\State Count = 0
\While{$tEST_1>tTOL$ and $Count< MaxRetry$}
\State $Count = Count +1$.
\State $k_{n+1} = \max\left\{ 0.9 k_n (\frac{tTOL}{tEST_1(n+1)})^{1/2}, 0.5 k_{n+1} \right\}$,
\State $k_{n+1} = \max\left\{ k_{n+1}, \Delta t_{\min}\right\}$.
\State REPEAT step
\State Break while loop if $k_{n+1}=\Delta t_{\min}$.
\EndWhile
\EndIf
\State {\bf Update the penalty parameter $\epsilon$ using the following tree:}
\State Compute divergence on each triangle 
\begin{equation*}
est_{\Delta}=\int_{\Delta} |\nabla \cdot u^{h}_{\epsilon,n+1}|^2\, dx.
\end{equation*}
\State Update $\epsilon_\Delta$:
\begin{equation*}
\begin{gathered}
\rho = \frac{LocTol_\Delta}{est_\Delta},\\
\epsilon_{\Delta,n+2} \leftarrow min (max(EpsMin, \rho\times \epsilon_{\Delta, n+1}), EpsMax).
\end{gathered}
\end{equation*}
Recover pressure $p$ if needed
\begin{equation*}
p_{\Delta, n+1} = -\frac{1}{\epsilon_{\Delta,n+1}} \nabla \cdot u^{h}_{\epsilon,n+1}.
\end{equation*}
\EndWhile
\end{algorithmic}
\end{algorithm}

\begin{figure}
    \centering
\includegraphics[width=0.85\textwidth]{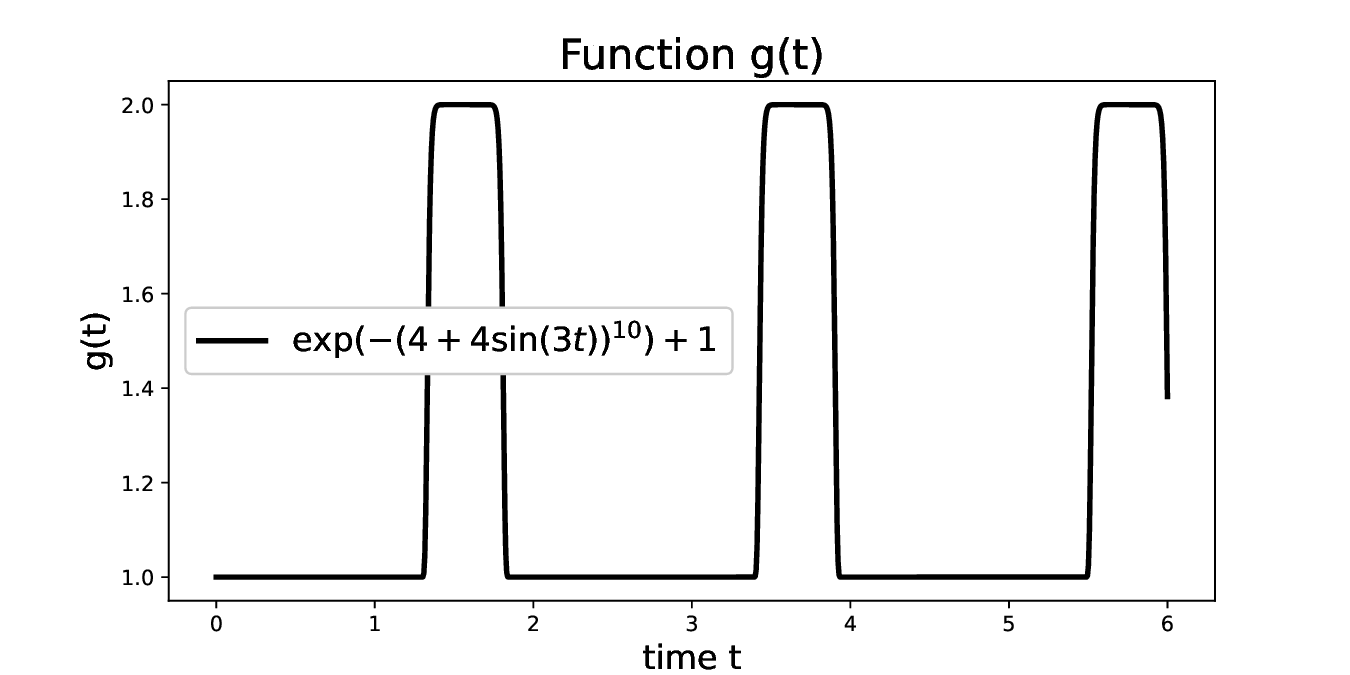}
    \caption{ Plot of $g(t)$.}
    \label{fig: sharp-g}
\end{figure}

\textit{Test 3: Body force with sharp transition regions without an exact solution} In this test, we create a body force with sharp transition regions in time within the domain $\Omega = (-1, 1) \times (-1, 1)$ in 2d. The body force $f$ is defined by 
\begin{equation*}f(x,y,t)= h(x,y,t) 
g(t),
\end{equation*}
where $h(x,y,t)$ is the body force from a test with an exact solution, and $g(t)$ helps to create the sharp transition regions in time. $h(x,y,t)$ is the body force for the following test with $\nu=1$ and $Re=1/\nu$ from \cite{decaria2017conservative},\cite{kean2023doubly}:
\begin{equation}\label{test-exacty-3}
\begin{gathered}
u(x,y,t) = \pi \sin(t) (\sin(2\pi y) \sin^2(\pi x), -\sin(2\pi x) \sin^2(\pi y))^T,\\
p(x,y,t) = \sin(t) \cos(\pi x) sin(\pi y).
\end{gathered}
\end{equation}
This $g(t)$ achieves sharp transition regions in time (see the $g(t)$ plot in Figure \ref{fig: sharp-g}): 
\begin{equation*}
g(t)= \exp(-\left(4+4\sin(3 t)\right)^{10})+1.
\end{equation*}
We generate unstructured meshes with GMSH \cite{geuzaine2009gmsh} using a target mesh size of $lc = 0.01$. We set an initial time step $\Delta t = 0.01$ and the final time $T = 6$. We impose the no-slip boundary conditions and set the initial condition $u_0=0.1$ in $\Omega$. $EpsMin=10^{-8}$, $EpsMax=10^{-1}$, and $TOL=10^{-5}$. The tolerance for time is set to $tTOL = 10^{-4}$, with a minimum tolerance of $tTOL_{\min} = tTOL/10$. We compare the locally adaptive penalty algorithm with adaptive time stepping to the locally adaptive penalty with a constant timestep. We apply the second order linear extrapolation of $u$ to the nonlinear term $b^*(u^*, u^{h}_{\epsilon, n+1}, v^h)$ for both adaptive time stepping and constant timestep, where 
\begin{equation*}
u^*=(1+\frac{k_{n+1}}{k_n})u^{h}_{\epsilon,n} -\frac{k_{n+1}}{k_n}u^{h}_{\epsilon,n-1}.
\end{equation*}
Note that in Algorithm \ref{alg: doubly-first-order}, we set a minimum time step of $\Delta t_{\min} = 0.001$, a maximum time step of $\Delta t_{\max} = 0.1$, and the maximum number of retries at a timestep $MaxRetry = 10$, which helps in controlling the time stepping behavior of the algorithm. 

\begin{figure}
    \centering
\includegraphics[width=0.85\linewidth]{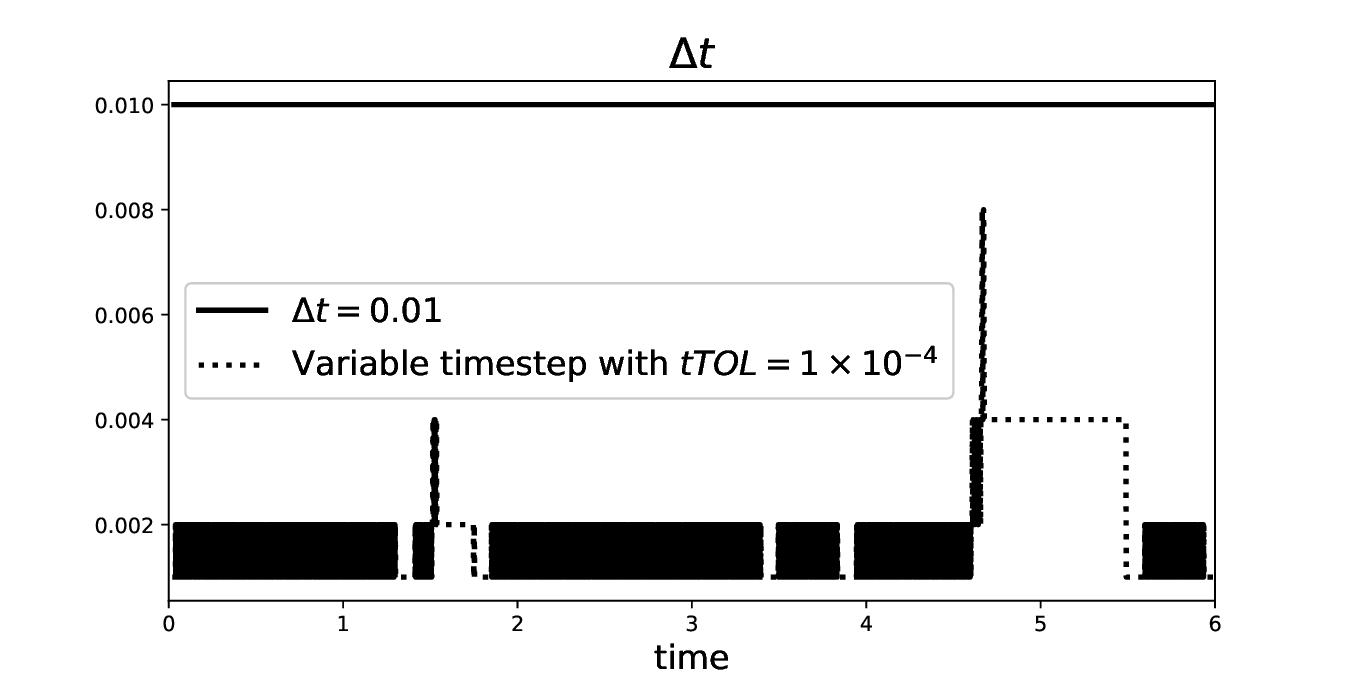}
    \caption{Timestep $\Delta t$ evolve in time.}
    \label{fig:time-step}
\end{figure}

\begin{figure}
    \centering
\includegraphics[width=0.85\linewidth]{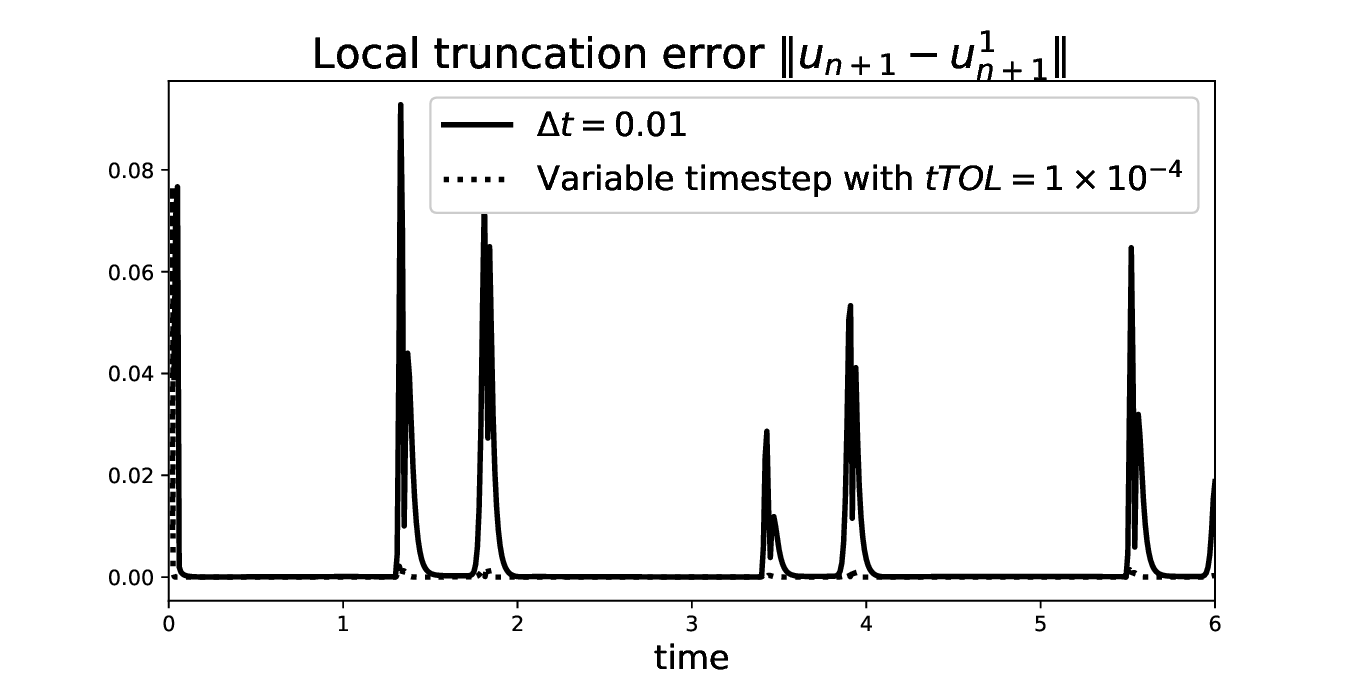}
    \caption{The local truncation error.}
    \label{fig:local-truncation-error}
\end{figure}

\begin{figure}
    \centering
\includegraphics[width=0.85\linewidth]{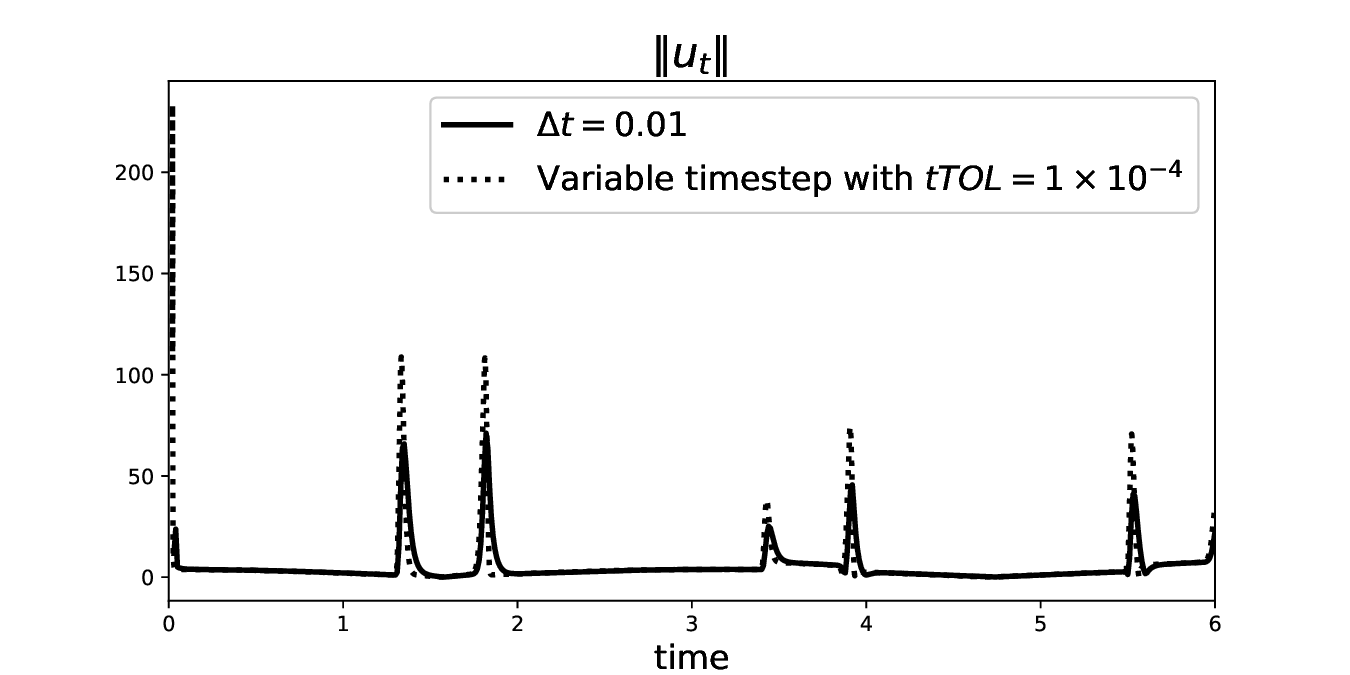}
    \caption{The time derivative of the fluid velocity $\|u_t\|$.}
    \label{fig:u_t}
\end{figure}

\begin{figure}
\centering
\includegraphics[width=0.85\linewidth]{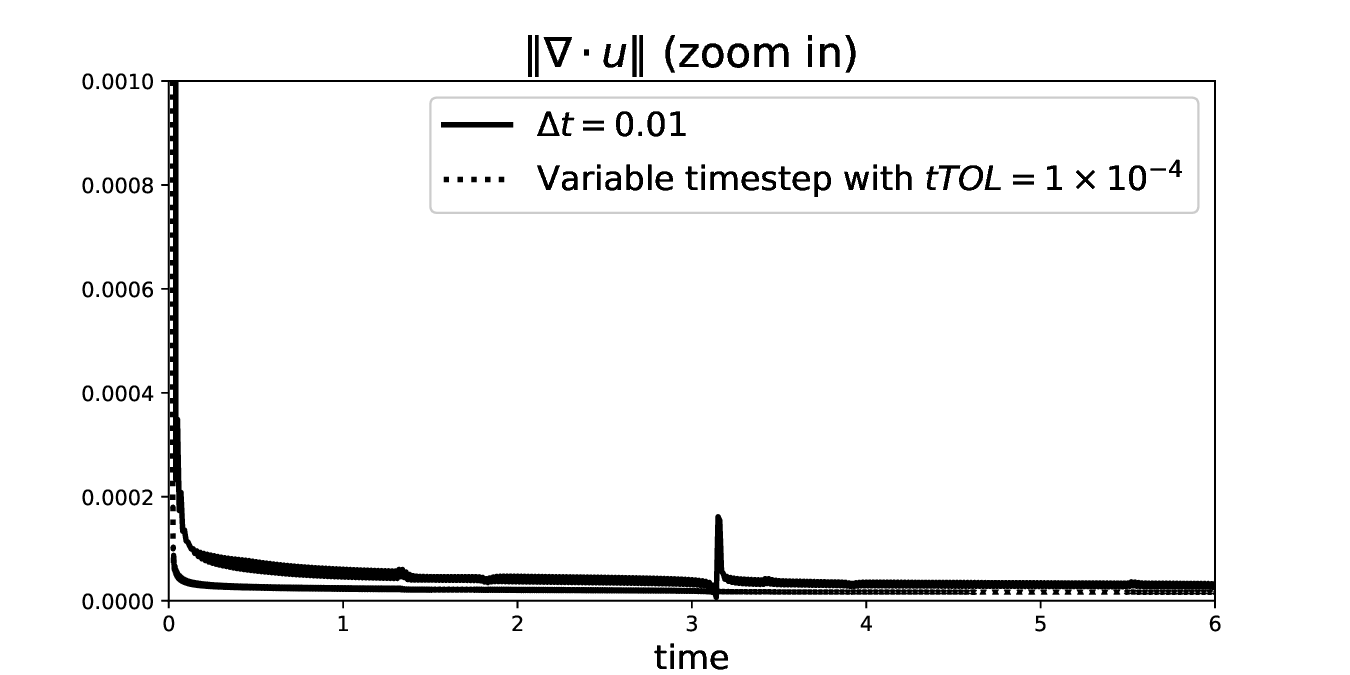}
    \caption{The divergence of velocity $\|\nabla \cdot u\|$.}
    \label{fig:div-test3}
\end{figure}

We plot the timestep $\Delta t$ over time in Figure \ref{fig:time-step}. The variable time stepping results in a smaller timestep to resolve the sharp transitions and ensure the local truncation error is bounded by $tTOL$. In Figure \ref{fig:local-truncation-error}, we can see smaller local truncation errors when using the locally adaptive penalty with adaptive time stepping, particularly around time intervals where the force changes rapidly.\\
The effectiveness of adaptive time stepping is highlighted by the larger $\|u_t\|$ values with better resolution of transitions near sharp transition regions, as shown in Figure \ref{fig:u_t}. The adaptive stepping method gives a larger fluid velocity derivative than the constant timestep method, which is more accurate.\\
Note that the penalty parameter is updated separately from the time stepping. Both constant time step and adaptive step methods give a good approximation of the incompressible flow in Figure \ref{fig:div-test3}. Variable time stepping method smooths out the divergence spike near $t=3.15$.\\
In summary, the locally adaptive method works well when combined with adaptive time stepping. When solving a flow problem with a sharp transition, the adaptive time-stepping method outperforms the constant time-stepping method. 
\section{Conclusions and open problems}\label{sec: conclusion}Penalty methods relax the incompressibility condition. This approach uncouples velocity and pressure and speeds up the calculation. Xie \cite{xie2022adaptive} developed an adaptive penalty scheme for the Stokes problem. Herein, we study the natural extension to the nonlinear, time-dependent NSE. We have proved the stability and derived an error estimate of the pointwise adaptive penalty on the NSE. We employed a modified green-Taylor vortex test on the convergence rate and observed $\|\nabla \cdot u\| = \mathcal{O}(TOL)$. Then we work on a complex flow problem and compare our algorithm with constant penalty $\epsilon= \Delta t$, and the non-penalized coupled system. Our algorithm is sufficient to control the divergence of velocity. We observe $\|\nabla \cdot u  \| \lesssim TOL$, Figure \ref{fig: div}. We combine the locally adaptive penalty method with the adaptive time stepping (see Algorithm \ref{alg: doubly-first-order}) and test it on a flow with sharp transition regions, Section \ref{sec: doubly-adaptive}. The doubly adaptive method (Algorithm \ref{alg: doubly-first-order}) improves the flow accuracy compared with the standard locally adaptive penalty method. One can combine the locally adaptive penalty method with adaptive time stepping with variable order methods in \cite{kean2023doubly} and extend this algorithm to the ensemble penalty scheme \cite{fang2023penalty}.

Pressure recovery for the locally adaptive $\epsilon$ penalty method is an open problem. Fairag \cite{fairag2002two}, Kean and Schneier \cite{kean2020error} gave some pressure recovery methods. Other open problems include an extension to coupled systems including a flow component and analysis of the effect of variability of $\epsilon$ on conditioning. 
 
\section*{Acknowledgments}
I thank Professor William Layton for his instruction, guidance, and support throughout this research. The author was partly supported by the NSF under grant DMS 2110379 and the University of Pittsburgh Center for Research Computing through the resources provided on the SMP cluster.

\bibliography{mybib}

\end{sloppypar}
\end{document}